\documentclass[12pt]{article}

\usepackage{amsfonts,amsmath}
\usepackage{latexsym}

\author{ K\'aroly J. B\"or\"oczky\footnote{Supported by
OTKA grants 068398 and 75016, and by the EU Marie Curie TOK
project DiscConvGeo and FP7 IEF grant GEOSUMSETS.}, Keith M. Ball}
\title{Stability of the Pr\'ekopa-Leindler inequality}

\newcommand{\proof}{\noindent{\it Proof: }}
\newcommand{\proofbox}{\mbox{ $\Box$}\\}
\newcommand{\R}{\mathbb{R}}

\newtheorem{lemma}{LEMMA}[section]
\newtheorem{theo}[lemma]{THEOREM}
\newtheorem{example}[lemma]{Example}
\newtheorem{coro}[lemma]{COROLLARY}
\newtheorem{prop}[lemma]{PROPOSITION}
\newtheorem{remark}[lemma]{REMARK}

\begin{document}

\maketitle

\begin{abstract}
We prove a stability version of the
Pr\'ekopa-Leindler inequality.
\end{abstract}

\section{The problem}

Our main theme is the Pr\'ekopa-Leindler inequality, due to
A. Pr\'ekopa \cite{Pre71} and L. Leindler \cite{Lei72}.
Soon after its proof, the inequality was generalized in
 A. Pr\'ekopa \cite{Pre73} and
\cite{Pre75}, C. Borell \cite{Bor75}, and in H.J. Brascamp, E.H. Lieb \cite{BrL76}.
Various applications are provided and surveyed in K.M. Ball \cite{Bal},
F. Barthe \cite{Bar}, and R.J. Gardner \cite{Gar02}.
The following multiplicative version from \cite{Bal}, is often more useful
and is more convenient for our purposes.

\begin{theo}[Pr\'ekopa-Leindler]
If $m,f,g$ are non-negative integrable functions on $\R$
satisfying $m(\frac{r+s}2)\geq \sqrt{f(r)g(s)}$ for
$r,s\in\mathbb{R}$, then
$$
\int_{\R}  m\geq \sqrt{\int_{\R}f \cdot \int_{\R}g}.
$$
\end{theo}

S. Dubuc \cite{Dub77} characterized the equality case  if the
integrals of
 $f,g,m$ above are positive. For this characterization,
we say that a non-negative real function $h$ on $\R$
is log-concave if for any $x,y\in\R$ and $\lambda\in(0,1)$,
we have
$$
h((1-\lambda) x+\lambda y)\geq h(x)^{1-\lambda}h(y)^{\lambda}.
$$
In other words, the support of $h$ is an interval, and $\log h$
is concave on the support. Now \cite{Dub77} proved that
equality holds in the Pr\'ekopa-Leindler inequality if and only if
there exist $a>0$, $b\in\R$ and a log-concave $h$ with positive
integral on $\R$ such that for a.e. $t\in\R$, we have
\begin{eqnarray*}
m(t)&=& h(t)\\
f(t)&=& a\cdot h(t+b)\\
g(t)&=& a^{-1}\cdot h(t-b).
\end{eqnarray*}
In addition for all $t\in R$, we have $m(t)\geq h(t)$, $f(t)\leq
a\cdot h(t+b)$ and $g(t)\leq a^{-1}\cdot h(t-b)$.

Our goal is to prove a stability version of the Pr\'ekopa-Leindler
inequality.

\begin{theo}
\label{PLstab}
There exists an positive absolute constant $c$ with
the following property.
If $m,f,g$ are non-negative integrable functions
with positive integrals on $\R$ such that
$m$ is log-concave,  $m(\frac{r+s}2)\geq \sqrt{f(r)g(s)}$ for
$r,s\in\mathbb{R}$,  and
$$
\int_{\R}  m\leq (1+\varepsilon) \sqrt{\int_{\R}f \cdot \int_{\R}g},
$$
for $\varepsilon>0$, then
there exist $a>0$, $b\in\R$ such that
\begin{eqnarray*}
\int_{\R}|f(t)-a\,m(t+b)|\,dt&\leq &
c\cdot\sqrt[3]{\varepsilon}|\ln \varepsilon|^{\frac43}\cdot
a\cdot\int_{\R}m(t)\,dt \\
\int_{\R}|g(t)-a^{-1}m(t-b)|\,dt&\leq &
c\cdot\sqrt[3]{\varepsilon}|\ln \varepsilon|^{\frac43}\cdot
a^{-1}\cdot\int_{\R}m(t)\,dt.
\end{eqnarray*}
\end{theo}

\begin{remark}
The statement also holds if the condition that $m$ is
log concave, is replaced by the condition
 that both $f$ and $g$ are log-concave. The reason
is that the function
$$
\tilde{m}(t)=\sup\{\sqrt{f(r)g(s)}:\,t=\mbox{$\frac{r+s}2$}\}
$$
is log-concave in this case.
\end{remark}

\begin{remark}
Most probably, the optimal
error estimate in Theorem~\ref{PLstab}
is of order $\varepsilon$. This cannot be proved
using the method of this note;
namely, by proving first
an estimate on the quadratic transportation distance.
\end{remark}

Let us summarize the main idea to prove
Theorem~\ref{PLstab}.
It can be assumed that $f$ and $g$ are
log-concave probability distributions
with zero mean (see Section~\ref{secmainproof}).
We establish the main properites
of log-concave functions in Section~\ref{logconcave},
and introduce the transportation map in Section~\ref{sectrans}.
After translating the condition
$\int_{\R}  m\leq (1+\varepsilon) \sqrt{\int_{\R}f \cdot \int_{\R}g}$
into an estimate for the
transportation map, we estimate the quadratic transportation distance
in Section~\ref{sectransest}.
Based on this, we estimate the $L_1$ distance of
$f$ and $g$ in Section~\ref{sectransL1},
which leads to the proof Theorem~\ref{PLstab}
in Section~\ref{secmainproof}. We note that
the upper bound
in Section~\ref{sectransL1} for the $L_1$ distance
of two log-concave probability distributions
in terms of the their quadratic  transportation distance
is close to being optimal.

\begin{remark}
It is not clear whether the condition in Theorem~\ref{PLstab}
that $m$ is log-concave is necessary
for there to be a stability estimate.
\end{remark}

\begin{remark}
Given $\alpha,\beta\in(0,1)$ with $\alpha+\beta=1$,
we also have the following version of the
Pr\'ekopa-Leindler inequality:
If $m,f,g$ are non-negative integrable functions on $\R$
satisfying $m(\alpha r+\beta s)\geq f(r)^\alpha g(s)^\beta$ for
$r,s\in\mathbb{R}$, then
$$
\int_{\R}  m\geq \left(\int_{\R}f\right)^\alpha
\left(\int_{\R}g\right)^\beta.
$$
The method of this note also yields the corresponding
stability estimate, only the $c$ in the new version
of Theorem~\ref{PLstab} depends on $\alpha$.
For this statement, the formula
$$
\frac{1+T'(x)}{2\sqrt{T'(x)}}= 1+
\frac{(1-T'(x))^2}{2\sqrt{T'(x)}(1+\sqrt{T'(x)})^2},
$$
used widely in this note
if $T'(x)$ is ``not too large'', should be replaced with
Koebe's  estimate
$$
\frac{\alpha+\beta T'(x)}{T'(x)^\beta}\geq 1+
\frac{\min\{\alpha,\beta\}(1-T'(x))^2}{T'(x)^\beta(1+\sqrt{T'(x)})^2}.
$$
In addition, if $T'(x)$ is ``large'', then one uses
$\frac{\alpha+\beta T'(x)}{T'(x)^\beta}>\beta T'(x)^\alpha$.
\end{remark}

\begin{remark}
The Pr\'ekopa-Leindler inequality
also holds in higher dimensions. One possible approach
to have a higher dimensional analogue of the stability
statement is to use  Theorem~\ref{PLstab}
and a suitable stability version of the injectivity
of the Radon transform on log-concave functions.
Here the difficulty is caused by the fact that the Radon transform
is notoriously instable even on the space of
smooth functions. Another possible approach is
to use recent stability version of the Brunn-Minkowski inequality
due to A. Figalli, F. Maggi, A. Pratelli \cite{FMP1} and \cite{FMP2}.
This approach has been successfully applied
in K.M. Ball, K.J. B\"or\"oczky \cite{BB10}.
\end{remark}

\section{Some elementary properties of log-concave
 probability distributions on $\R$}
\label{logconcave}

Let $h$ be a log-concave probability distribution on $\R$. In this
section we discuss various useful elementary properties of $h$.
Many of these properties are implicit or explicit in
\cite{Bal88}.

First, assuming $h(t_0)=a\cdot b^{t_0}$ for $a,b>0$, and
$t_1<t_0<t_2$, we have
\begin{equation}
\label{logconcavity}
\begin{array}{l}
\mbox{if $h(t_1)\geq a\cdot b^{t_1}$, then
 $h(t_2)\leq a\cdot b^{t_2}$,}\\
\mbox{if $h(t_2)\geq a\cdot b^{t_2}$, then
 $h(t_1)\leq a\cdot b^{t_1}$}.
\end{array}
\end{equation}
Next we write $w_h$ and $\mu_h$ to denote the median and mean
of $h$;
namely,
$$
\int_{-\infty}^{w_h}h=\int_{w_h}^{\infty}h=\frac12
\mbox{ \ and \ }\mu_h=\int_{\R}xh(x)\,dx.
$$

\begin{prop}
\label{moment}
If $f$ and $g$ are positive, and $\theta$ is an
increasing function on $(a,b)$, and there exists $c\in(a,b)$ such
that $f(t)\leq g(t)$ if $t\in(a,c)$, and $f(t)\geq g(t)$ if
$t\in(c,b)$, and $\int_a^bg(t)\,dt=\int_a^bf(t)\,dt$ then
$$
\int_a^b\theta(t)g(t)\,dt\leq \int_a^b\theta(t)f(t)\,dt.
$$
\end{prop}
\proof We may assume that $g(t)>0$ if $t\in(a,c)$, and $f(t)>0$ if
$t\in(c,b)$. Let $(a_0,b)$ and $(a,b_0)$ be the support of $f$ and
$g$, respectively, where $a_0\in[a,c]$ and $b_0\in[c,b]$. Let
$T:(a_0,b)\to (a,b_0)$ be the transportation map defined by
$$
\int_{a_0}^xf(t)\,dt= \int_{a}^{T(x)}g(t)\,dt.
$$
The conditions yield that $T$ is monotone increasing, bijective,
continuous, $T(x)\leq x$ for $x\in(a_0,b)$, and for a.e. $x\in
(a_0,b)$, we have
$$
f(x)=g(T(x)) T'(x).
$$
Therefore
\begin{eqnarray*}
\int_a^b\theta(t)g(t)\,dt&=&\int_{a_0}^b\theta(T(s))g(T(s))T'(s)\,ds=
\int_{a_0}^b\theta(T(s))f(s)\,ds\\
 &\leq&
\int_a^b\theta(s)f(s)\,ds. \proofbox
\end{eqnarray*}

\begin{prop}
\label{hwestimate}
If $h$ is a log-concave probability distribution on
$\R$ then for $w=w_h$ and $\mu=\mu_h$, we have
\begin{description}
\item{(i)}
$h(w)\cdot |w-\mu|\leq \ln\sqrt{e/2}$.
\item{(ii)}
$ h(w)\cdot e^{-2h(w)|x-w|}\leq h(x)\leq h(w)\cdot e^{2h(w)|x-w|}
\mbox{ \ \ if $|x-w|\leq\frac{\ln 2}{2h(w)}$}$.
\item{(iii)}
$h(x)\leq 2h(w)$ for $x\in \R$.
\item{(iv)} If $x>w$ then
$\int_x^{\infty}h\leq\frac{h(x)}{2h(w)}$.
\item{(v)} If $x>w$ and $\int_x^{\infty}h=\nu>0$, then
\begin{eqnarray*}
\int_x^{\infty}(t-w) h(t)\,dt&\leq& \frac{\nu}{4h(w)}\cdot(1-\ln 2\nu)\\
\int_x^{\infty}(t-w)^2 h(t)\,dt&\leq& \frac{\nu}{8h(w)^2}\cdot
[(\ln 2\nu)^2-2\ln 2\nu+2].
\end{eqnarray*}
\end{description}
\end{prop}
{\bf Remark } All estimates are optimal.\\
\proof We may assume that $w_h=0$, and $h(w)=\frac12$. It is
natural to compare $h$ near $0$ to the probability distribution
$$
\varphi(x)=\left\{
\begin{array}{lcl}
\frac12\cdot e^{-x}&\mbox{ }&
\mbox{if $x\geq -\ln 2$}\\[0.5ex]
0&\mbox{ }& \mbox{if $x< -\ln 2$},
\end{array}
\right.
$$
which satisfies $w_\varphi=0$, and $\varphi(0)=h(0)$.
 Since
$\int_{-\infty}^{0}h=\int_{-\infty}^{0}\varphi$, we have $h(x)\geq
\varphi(x)$ for some $x>0$.
 It follows from (\ref{logconcavity})
that there exists some $v>0$ such that
\begin{equation}
\label{compareh}
\begin{array}{rclcl}
h(x)&\geq & \varphi(x)&\mbox{ }&\mbox{provided $x\in[0,v]$}\\
h(x)&\leq & \varphi(x) &\mbox{ }&\mbox{provided $x\geq v$ or
$x\in[-\ln 2,0]$}.
\end{array}
\end{equation}
In particular
$\int_{-\infty}^{0}h=\int_{-\infty}^{0}\varphi$
and $\int_{0}^{\infty}h=\int_{0}^{\infty}\varphi$ yield
\begin{eqnarray*}
-\ln \frac{e}2&=&\int_{-\infty}^{0}x\varphi(x)\,dx
+\int_{0}^{\infty}x\varphi(x)\,dx\\
&\leq &\int_{-\infty}^{0}xh(x)\,dx
+\int_{0}^{\infty}xh(x)\,dx =\mu.
\end{eqnarray*}
Comparing $h$ to $\varphi(-x)$ shows that
$\mu\leq \ln \frac{e}2$, and in turn, we deduce (i).

Turning to (ii), the upper bound directly follows from
(\ref{compareh}). To prove the lower bound, we may assume that
$x>0$. According to (\ref{compareh}), it is enough to check the
case $x=\ln 2$. Therefore we suppose that
$$
h(\ln 2)<1/4,
$$
and seek a contradiction. Since $h$ is log-concave, there exists
some $a\in\R$ such that
$$
h(x)< \mbox{$\frac1{4}\,e^{-a(x-\ln 2)}$ \ for $x\in\R$}.
$$
Here $h(0)=\frac12$ yields that $a\geq 1$.

We observe that $\frac1{4}\,e^{a(x_0-\ln 2)}=\frac12\,e^{x_0}$ for
$x_0=\frac{a-1}{a+1}\,\ln 2$, and applying the analogue of
(\ref{compareh}) to $\varphi(-x)$, we obtain that $h(x)\leq
\frac12\,e^{x}$ for $x\in [0,x_0]$.
 In particular
$$
\int_0^{\infty}h<\int_0^{x_0}\frac12\, e^{x}dx+
\int_{x_0}^{\infty}\frac1{4}\,e^{-a(x-\ln 2)}dx=
\left(\frac1a+1\right)2^{-\frac2{a+1}}-\frac12.
$$
Differentiation shows that the last expression is first
decreasing, and after increasing in $a\geq 1$. Since the value of
this last expression is $\frac12$ both at $a=1$ and at $a=\infty$,
we deduce that $\int_0^{\infty}h<\frac12$. This is absurd,
therefore we have proved (ii).

To prove (iii), we may assume $x>0$ and $h(x)\geq 1$. Since
$h(t)\geq \frac12\,e^{\frac{t}{x}\, \ln 2h(x)}$ for $t\in[0,x]$,
we have
$$
\frac12\geq \int_0^xh\geq \int_0^x
\frac12\,e^{\frac{t}{x}\, \ln 2h(x)}dt =
\frac{x(2h(x)-1)}{2\ln 2h(x)}.
$$
Now (ii) and $h(x)\geq 1$ yield that $x\geq \ln 2$. As
$\frac{s-1}{s}>\frac1{\ln 2}$ for $s>2$, we conclude $h(x)\leq 1$.

To prove (iv), we may assume that
$h(x)<h(w)$. Let $x_0=-\ln 2h(x)$, and hence
$h(x)=\frac12\,e^{-x_0}$. If $h(x)\geq h(v)$, then
(\ref{compareh}) yields
 $\int_0^{x}h(t)\,dt\geq \int_0^{x_0}\frac12\,e^{-t}\,dt$,
 and hence
$$
\int_x^{\infty}h(t)\,dt\leq
\int_{x_0}^{\infty}\frac12\,e^{-t}\,dt=h(x).
$$

If $h(x)<h(v)$ then $x_0\geq x$. We may assume
 that $x_0>x$, and hence $h(x)<\frac12\,e^{-x}$.
 We choose $a>0$ such that
 $$
\int_x^{\infty}h(t)\,dt= \int_{x}^{\infty}h(x)\,e^{-a(t-x)}\,dt,
 $$
and consider the function
$$
\tilde{h}(t)=\left\{
\begin{array}{lcl}
h(t)&\mbox{ }&
\mbox{if $t\leq x$}\\[0.5ex]
h(x)\,e^{-a(t-x)}&\mbox{ }& \mbox{if $t\geq x$},
\end{array}
\right.
$$
If follows by the choice of $a$ that $h(t)\geq \tilde{h}(t)$
 for some $t>x$. We deduce that $\tilde{h}$ is also log-concave,
 and hence $\tilde{h}(t)\leq \frac12\,e^{-t}$ for $t\geq v$.
 Therefore $a\geq 1$, and we conclude
 that
$$
\int_x^{\infty}h(t)\,dt=\int_x^{\infty}h(x)\,e^{-a(t-x)}\,dt=h(x)/a\leq
h(x).
$$

Finally, we prove (v). Let $x_1=-\ln 2\nu$, which satisfies that
$\int_x^{\infty}h(t)\,dt=\int_{x_1}^{\infty}\frac12\,e^{-t}\,dt$.
It follows from  (\ref{compareh}) that $x_1\geq x$.
We define two functions $f$ and $g$ on $[x,\infty)$.
Let $f(t)=\frac12\,e^{-t}$ if $t\geq x_1$,
and let $f(t)=0$ if $t\in[x,x_1)$. In addition
let $g=h|_{[x,\infty)}$. These two functions satisfy
the conditions in Proposition~\ref{moment},
therefore for $\alpha\geq 0$, we have
$$
\int_x^{\infty}t^{\alpha}h(t)\,dt=
\int_x^{\infty}t^{\alpha} g(t)\,dt\leq
\int_x^{\infty}t^{\alpha}f(t)\,dt=
\int_{x_1}^{\infty} \frac{t^{\alpha}e^{-t}}2\,dt.
$$
Evaluating the last integral for $\alpha=1,2$
yields (v). \proofbox


Next we discuss various consequences of
Proposition~\ref{hwestimate}.

\begin{coro}
\label{hwest}
Let $h$ be a log-concave probability density function on $\R$,
and let $\int_x^{\infty}h=\nu\in(0,\frac12]$. Then
\begin{description}
\item{(i)} $h(x)\cdot e^{-\frac{h(x)|t-x|}{\nu}}
\leq h(t)\leq h(x)\cdot e^{\frac{h(x)|t-x|}{\nu}}$
 \ \ if $|t-x|\leq \frac{\nu\ln 2}{h(x)}$;
\item{(ii)} If $\nu\in(0,\frac16)$,
$w=w_h$ and $\mu=\mu_h$, then
\begin{eqnarray*}
\int_x^{\infty}|t-\mu| h(t)\,dt&\leq& \frac{\nu}{2h(w)}\cdot|\ln \nu|\\
\int_x^{\infty}|t-\mu|^2 h(t)\,dt&\leq& \frac{5\nu}{4h(w)^2}\cdot
(\ln \nu)^2.
\end{eqnarray*}
\end{description}
\end{coro}
{\bf Remark } The order of all estimates is optimal, as it is shown
by the example of $h(t)=e^{-|t|}/2$.\\
\proof To prove (i), let $|t-x|\leq \frac{\nu\ln 2}{h(x)}$.
There exists
a unique $\lambda\in\R$, such that for the function
$$
\tilde{h}(t)=\left\{
\begin{array}{ll}
h(t) &  \mbox{ \ \ if $t\geq x$;}\\
\min\{h(t),h(x)\cdot e^{\lambda(t-x)}\} &  \mbox{ \ \ if $t\leq x$.}
\end{array}
\right.,
$$
we have $\int_{-\infty}^x\tilde{h}=\nu$.
We note that $\tilde{h}$
is log-concave, and $\lambda\geq\frac{-h(x)}{\nu}$.
In particular $\frac1{2\nu}\,\tilde{h}$
is a log-concave probability distribution whose
median is $x$, and hence Proposition~\ref{hwestimate} (ii)
yields $h(t)\geq \tilde{h}(t)\geq h(x)\cdot e^{\frac{-h(x)|t-x|}{\nu}}$.
Since for $s=2x-t$, we have
$h(s)\geq h(x)\cdot e^{\frac{-h(x)|s-x|}{\nu}}$,
we conclude (i) by (\ref{logconcavity}).

For (ii), we may assume that $h(w)=\frac12$, and hence
Proposition~\ref{hwestimate} (i) yields
that $|w-\mu|\leq \ln\frac{e}2$. Since
$\ln 2\nu\leq -1$, we deduce by Proposition~\ref{hwestimate} (ii)
that
\begin{eqnarray*}
 \int_x^{\infty}|t-\mu| h(t)\,dt&\leq&
\int_x^{\infty}[|t-w|+|w-\mu|] h(t)\,dt\\
&\leq& \mbox{$
\nu\cdot(-\ln 2\nu)+
\nu\cdot\ln\frac{e}2<\nu\cdot|\ln \nu|$}.
\end{eqnarray*}
In addition
\begin{eqnarray*}
 \int_x^{\infty}(t-\mu)^2 h(t)\,dt&\leq&
\int_x^{\infty}2[(t-w)^2+(w-\mu)^2] h(t)\,dt\\
&\leq& \mbox{$
\nu\cdot 5 (\ln 2\nu)^2+
\nu\cdot 2(\ln\frac{e}2)^2<5\nu\cdot(\ln \nu)^2$}.\proofbox
\end{eqnarray*}

\section{The transportation map for log-concave
probability distributions, and
the Pr\'ekopa-Leindler inequality}
\label{sectrans}

Let $f$ and $g$ be log-concave probability distributions on $\R$,
and let $I_f$ and $I_g$ denote the open intervals that are the
supports of $f$ and $g$, respectively. We define the
transportation map $T:\,I_f\to I_g$ by the identity
\begin{equation}
\label{transdef}
\int_{-\infty}^xf(t)\,dt= \int_{-\infty}^{T(x)}g(t)\,dt.
\end{equation}
 In particular $T$ is monotone increasing,
bijective, and continuous on $I_f$, and for a.e. $x\in I_f$, we
have
\begin{equation}
\label{transprop}
f(x)=g(T(x)) T'(x).
\end{equation}

\noindent{\bf Remark } Using (\ref{transdef}),
the transportation map $T:\R\to\R$
can be defined for any two probability
distributions $f$ and $g$, and $T$ is
naturally monotone increasing.
It  is easy to see that (\ref{transprop}) holds
if there exists a set $A\subset\R$ of zero measure
such that both $f$ and $g$ are continous on
$\R\backslash A$.
Unfortunately
(\ref{transprop}) does not hold
in general. Let $B\subset\R$ be
such a set that the density of each point
of $B$ is strictly between $0$ and $1$,
and let $f$ be a probability distribution
that is zero on $\R\backslash B$,
and positive and continuous on $B$.
If say $g(x)=\frac12\,e^{|x|}$, then (\ref{transprop})
never holds. \\

For an integrable function $m$ on $\R$ satisfying
$m(\frac{r+s}2)\geq \sqrt{f(r)g(s)}$ for $r,s\in\mathbb{R}$, one
proof of the  Pr\'ekopa-Leindler inequality runs as follows:
\begin{eqnarray*}
1&=&\int_{\R} f=\int_{I_f} \sqrt{f(x)}\cdot
\sqrt{g(T(x))T'(x)}\,dx \\
&\leq&
\int_{I_f} m\left(\frac{x+T(x)}2\right)\sqrt{T'(x)}\,dx\\
&\leq&
\int_{I_f} m\left(\frac{x+T(x)}2\right)\cdot\frac{1+T'(x)}2\,dx\\
&=&\int_{\frac12(I_f+I_g)} m(x)\,dx\leq \int_{\R}m.
\end{eqnarray*}

The basic fact that we will exploit is this.
If we know that $\int_{\R}  m\leq 1+\varepsilon$ then
\begin{eqnarray}
\nonumber
\varepsilon&\geq &
\int_{I_f}m\left(\frac{x+T(x)}2\right)\cdot
\left(\frac{1+T'(x)}2-\sqrt{T'(x)}\right)dx\\
\nonumber
&\geq &\int_{I_f} \sqrt{f(x)}\cdot \sqrt{g(T(x))T'(x)}
\left(\frac{1+T'(x)}{2\sqrt{T'(x)}}-1\right)dx\\
\label{stabcond} &=&\int_{I_f}
f(x)\cdot\frac{(1-\sqrt{T'(x)})^2}{2\sqrt{T'(x)}}\,dx .
\end{eqnarray}

As long as $T'$ is not too large, the integrand is at least
about $f(x)(1-T'(x))^2$
and using a Poincar\'{e} inequality for the density $f$ we can
bound the integral
of this expression from below by the transportation cost $\int f(x)(x-T(x))^2$.
The main technical issue is to handle the places where $T'$ is large.

\section{The quadratic transportation distance}
\label{sectransest}

Let $f$ and $g$ be log-concave probability distributions on $\R$
with zero mean; namely,
$$
0=\int_{\R}xf(x)\,dx=\int_{\R}yg(y)\,dy.
$$
In this section we show that (\ref{stabcond}) yields
  an upper bound for the quadratic transportation distance
$$
\int_{I_f}f(x)(T(x)-x)^2dx
$$
of $f$ and $g$.

\begin{lemma}
\label{transdist}
If $f$ and $g$ are log-concave probability distributions on $\R$
with zero mean, and (\ref{stabcond}) holds for
$\varepsilon\in(0,\frac1{48})$, then
$$
\int_{I_f}f(x)(T(x)-x)^2dx \leq 2^{20} f(w_f)^{-2}\cdot
\varepsilon|\ln \varepsilon|^2.
$$
\end{lemma}
{\bf Remark } The optimal power of $\varepsilon$
is most probably $\varepsilon^2$ in Lemma~\ref{transdist}
(compare Example~\ref{exa2}).
For a possible proof of an improved
estimate, we should improve
on (\ref{bobkov}) if $R(x)=T(x)-x$ where $T$
is the transportation map  for another log-concave
probability distribution. One may possibly use that $T(x)-x$ is of
at most logarithmic order.\\
\proof The main tool in the proof of Lemma~\ref{transdist} is the Poincar\'{e}
inequality for log-concave measures which can be found in
 (1.3) and (4.2) of S.G. Bobkov \cite{Bob99}. If $h$
is a log-concave probability distribution on $\R$, and the
function $R$ on $\R$ is locally Lipschitz with
expectation $\mu=\int_{\R}h(x)
R(x)\,dx$, then
\begin{equation}
\label{bobkov}
\int_{\R}h(x) (R(x)-\mu)^2\,dx=
\int_{\R}h(x) R(x)^2\,dx-\mu^2\leq
h(w_h)^{-2}\cdot\int_{\R}h(x) R'(x)^2\,dx.
\end{equation}

We may assume that $g(w_g)\leq f(w_f)$, and
$f(w_f)=\frac12$. Let $T$ be the transportation map
from $f$ to $g$, and let $S$ be its inverse,
thus for a.e. $x\in I_f$ and $y\in I_g$, we have
\begin{equation}
\label{Ttrans}
 f(x)=g(T(x)) T'(x)\mbox{ \ and \ }
g(y)=f(S(y))S'(y).
\end{equation}

Suppose that for some
$x\in\R$ with $\int_x^\infty f=\nu\in(0,\frac12]$, we have
$g(T(x))\leq \frac1{16}\,f(x)$.
If $x\leq t\leq x+\frac{\nu\ln 2}{f(x)}$
then Corollary~\ref{hwest} (i) yields
$f(t)\geq f(x)\cdot e^{-\frac{f(x)(t-x)}{\nu}}\geq \frac12\,f(x)$.
On the other hand, the log-concavity of $g$ and
Proposition~\ref{hwestimate} (iii) yield that
if $x\leq t< x+\frac{\nu\ln 2}{f(x)}$, then
$g(t)<2g(x)\leq \frac1{4}\, f(t)$.
In particular $T'(t)>4$ by (\ref{Ttrans}),
and hence (compare (\ref{stabcond}))
$$
\varepsilon\geq
\int_{\R}\frac{(1-\sqrt{T'(t)})^2}{2\sqrt{T'(t)}}\, f(t)\,dt
>\int_x^{x+\frac{\nu\ln 2}{f(x)}}
\frac{f(x)}4\cdot e^{-\frac{f(x)(t-x)}{\nu}}\,dt
=\frac{\nu}8.
$$
Similar argument for $f(-x)$ and $g(-x)$
shows that if $\int_{-\infty}^x f=\nu$ and
$g(T(x))\leq \frac1{16}\,f(x)$ then $\nu<8\varepsilon$.

We define $x_1,x_2,y_1,y_2$ by
$$
\int_{-\infty}^{x_1} f=\int_{x_2}^\infty f=
\int_{-\infty}^{y_1} g=\int_{y_2}^\infty g=8\varepsilon<\frac16.
$$
The argument above yields that if
$x\in (x_1,x_2)$, then $T'(x)\leq 16$
and $g(T(x))\geq \frac1{16}\,f(x)$, and hence $g(w_g)\geq \frac1{32}$.
As the means of $f$ and $g$ are zero,
we deduce by Corollary~\ref{hwest} (ii)
and (\ref{Ttrans}) that
\begin{eqnarray}
\label{xf}
\int_{\R\backslash[x_1,x_2]}|x|f(x)\,dx&\leq&2^4\varepsilon|\ln\varepsilon|;\\
\label{Tf}
\int_{\R\backslash[x_1,x_2]}|T(x)|f(x)\,dx&=&
\int_{\R\backslash[y_1,y_2]}|y|g(y)\,dy
\leq 2^8\varepsilon|\ln\varepsilon|;\\
\label{x2f}
\int_{\R\backslash[x_1,x_2]}x^2f(x)\,dx&\leq&2^7\varepsilon(\ln\varepsilon)^2;\\
\label{T2f}
\int_{\R\backslash[x_1,x_2]}T(x)^2f(x)\,dx&=&
\int_{\R\backslash[y_1,y_2]}y^2g(y)\,dy
\leq 2^{15}\varepsilon(\ln\varepsilon)^2.
\end{eqnarray}
Since $(T(x)-x)^2\leq 2[T(x)^2+x^2]$, we have
\begin{equation}
\label{(T-x)2tail}
\int_{\R\backslash[x_1,x_2]}(T(x)-x)^2f(x)\,dx\leq
 2^{17}\varepsilon(\ln\varepsilon)^2.
\end{equation}
Next we consider the log-concave probability distribution
$$
\tilde{f}(t)=\left\{
\begin{array}{lcl}
(1-16\varepsilon)^{-1}f(t)&\mbox{ }&
\mbox{if $t\in [x_1,x_2]$}\\[0.5ex]
0&\mbox{ }& \mbox{if $t\in \R\backslash[x_1,x_2]$}.
\end{array}
\right.
$$
To estimate $\tilde{f}(w_{\tilde{f}})$,
we define $z_1=w_f-\ln 2$, and $z_2=w_f+\ln 2$.
Since $f(w_f)=\frac12$, Proposition~\ref{hwestimate} (ii)
applied to $f$ yields
$$
\int_{\R\backslash[z_1,z_2]}\tilde{f}(x)\,dx\leq
(1-16\varepsilon)^{-1}\left(1-16\varepsilon-
\int_{z_1}^{z_2}\frac{e^{-|x-w_f|}}2\,dx\right)<\frac12.
$$
It follows that $|w_{\tilde{f}}-w_f|<\ln 2$,
and hence we deduce again by Proposition~\ref{hwestimate} (ii)
that
$$
\tilde{f}(w_{\tilde{f}})> \frac14.
$$
For the expectation
$$
\mu=\int_\R(T(x)-x)\tilde{f}(x)\,dx,
$$
we have the estimate
$$
|\mu|=(1-16\varepsilon)^{-1}
\left|\int_{\R\backslash[x_1,x_2]}(T(x)-x)f(x)\,dx\right|\leq
2^{10}\varepsilon|\ln\varepsilon|.
$$
If
$x\in (x_1,x_2)$, then $T'(x)\leq 16$, thus the expression
in (\ref{stabcond}) satisfies
$$
\frac{(1-\sqrt{T'(x)})^2}{2\sqrt{T'(x)}}=
\frac{(T'(x)-1)^2}{2(1+\sqrt{T'(x)})^2\sqrt{T'(x)}}\geq
\frac{(T'(x)-1)^2}{200}>2^{-8}(T'(x)-1)^2.
$$
We deduce using (\ref{bobkov}) and (\ref{stabcond}) that
\begin{eqnarray}
\nonumber
\int_{[x_1,x_2]}(T(x)-x)^2f(x)\,dx&\leq&
\int_\R(T(x)-x)^2\tilde{f}(x)\,dx\\
\nonumber
&\leq& \mu^2+\tilde{f}(w_{\tilde{f}})^{-2}
\int_\R(T'(x)-1)^2\tilde{f}(x)\,dx\\
\nonumber
&\leq& 2^{20}\varepsilon^2|\ln\varepsilon|^2
+2^{13}\int_{x_1}^{x_2}\frac{(1-\sqrt{T'(x)})^2}{2\sqrt{T'(x)}}
\,f(x)\,dx\\
\label{middletrans}
&\leq&
2^{20}\varepsilon^2|\ln\varepsilon|^2+2^{13}\varepsilon.
\end{eqnarray}
However
(\ref{x2f}) and (\ref{T2f}) imply
$$
\int_{\R\backslash[x_1,x_2]}(T(x)-x)^2f(x)\,dx\leq
\int_{\R\backslash[x_1,x_2]}2(T(x)^2+x^2)f(x)\,dx\leq
 2^{17}\varepsilon(\ln\varepsilon)^2.
$$
Combining this estimate with (\ref{middletrans}),
completes the proof of Lemma~\ref{transdist}.
\proofbox

\section{The $L_1$ and
quadratic transportation distances}
\label{sectransL1}

Our goal is to estimate the $L_1$ distance of two log-concave
probability distributions $f$ and $g$ in terms of their quadratic
transportation distance. In this section, $T$ always denotes the
transportation map $T:\,I_f\to I_g$ satisfying
$$
\int_{-\infty}^xf(t)\,dt= \int_{-\infty}^{T(x)}g(t)\,dt.
$$
We prepare our estimate Theorem~\ref{transL1}
by Propositions~\ref{transL11} and \ref{transL10}.

When we write $A\ll B$ for expressions $A$ and $B$, then we mean
that $|A|\leq c\cdot B$ where $c>0$ is an absolute constant, and
hence is independent from all the quantities occurring  in $A$ and
$B$. In addition $A\approx B$ means that $A\ll B$ and $B\ll A$.

\begin{prop}
\label{transL11}
Let $f$ and $g$ be log-concave probability
distributions on $\R$ satisfying
$\int_{-\infty}^{z}f\geq \nu$ and $\int_{z}^\infty f\geq \nu$
for $\nu>0$ and $z\in\R$. If either $\int_{-\infty}^{z}g\leq \nu/2$ or
$\int_{z}^\infty g\leq \nu/2$, then
$$
\int_{z-\frac{\nu}{f(z)}}^{z+\frac{\nu}{f(z)}}(T(x)-x)^2f(x)\,dx\gg
\frac{\nu^3}{f(z)^2}.
$$
\end{prop}
\proof We may assume that $\int_{z}^\infty g\leq \nu/2$.
It follows from Corollary~\ref{hwest} (i)
that if $x\leq z+\frac{\nu\ln \frac{3}2}{f(z)}$
then $\int_z^x f\leq\nu/2$, and hence $T(x)\leq z$.
Therefore
$$
\int_{z+\frac{\nu\ln \frac{5}4}{f(z)}}^
{z+\frac{\nu\ln\frac{3}2}{f(z)}}(T(x)-x)^2f(x)\,dx\gg
\int_{z+\frac{\nu\ln \frac{5}4}{f(z)}}^
{z+\frac{\nu\ln\frac{3}2}{f(z)}}
\left(\frac{\nu\ln \frac{5}4}{f(z)}\right)^2\frac{f(z)}2\,dx
\gg \frac{\nu^3}{f(z)^2}. \mbox{ \ }\proofbox
$$

\begin{prop}
\label{transL10}
Let $f$ and $g$ be log-concave probability
distributions on $\R$ satisfying
$\int_{-\infty}^{z}f\geq \nu$ and $\int_{z}^\infty f\geq \nu$,
moreover $\int_{-\infty}^{z}g\geq \nu/2$ and $\int_{z}^\infty g\geq \nu/2$
for $\nu>0$ and $z\in\R$.
 If $g(z)\neq f(z)$ and
$\Delta=\frac{\nu\ln 2}{3f(z)}\cdot
\min\{|\ln\frac{g(z)}{f(z)}|,3\}$, then
$$
\int_{z-\Delta}^{z+\Delta}(T(x)-x)^2f(x)\,dx\gg\frac{\nu^3}{f(z)^2}\cdot
\min\left\{\left|\ln\frac{g(z)}{f(z)}\right|,3\right\}^4.
$$
\end{prop}
{\bf Remark} If
in addition $e^{-3}f(z)\leq g(z)\leq e^3 f(z)$,
then the arguments in Cases 2 and 3 show that
the interval $[z-\Delta,z+\Delta]$ of
integration can be replaced by
$[z-\frac{\Delta}{150},z+\frac{\Delta}{150}]$,
and if $x\in[z-\frac{\Delta}{150},z+\frac{\Delta}{150}]$, then
$$
\mbox{$\frac13\,|\ln\frac{g(z)}{f(z)}|
\leq|\ln\frac{g(x)}{f(x)}|\leq \frac53\,|\ln\frac{g(z)}{f(z)}|$}.
$$
\proof According to Corollary~\ref{hwest} (i),
if $z-\Delta\leq x\leq z+\Delta$,
then
\begin{equation}
\label{f(x)nu}
f(z)/2\leq f(z)\cdot e^{\frac{-f(z)|x-z|}{\nu}}\leq f(x)\leq
f(z)\cdot e^{\frac{f(z)|x-z|}{\nu}}\leq 2f(z).
\end{equation}
Similarly if $z-\frac{\nu\ln 2}{2g(z)}\leq x\leq z+\frac{\nu\ln 2}{2g(z)}$,
then
\begin{equation}
\label{g(x)nu}
g(z)/2\leq g(z)\cdot e^{\frac{-2g(z)|x-z|}{\nu}}\leq g(x)\leq
g(z)\cdot e^{\frac{2g(z)|x-z|}{\nu}}\leq 2g(z).
\end{equation}
We may assume
$$
T(z)\leq z.
$$
For the rest of the
argument, we distinguish four cases.\\

\noindent{\bf Case 1 } $g(z)\geq e^3 f(z)$.\\
In this case, $\Delta=\frac{\nu\ln 2}{f(z)}$.
We note that,
\begin{equation}
\label{numbers}
\frac{\ln 2}{2\cdot e^{3}}<\frac{\ln 2}{10}<\frac{3\ln 2}{10}
<\ln\frac{5}4.
\end{equation}
Since $\frac{\nu\ln 2}{2g(z)}<\frac{\Delta}{10}$,
(\ref{g(x)nu}) yields that if $x\geq z+\frac{\Delta}{10}$, then
\begin{equation}
\label{intgbig}
\int_{z}^{x}g>\frac{\nu}4.
\end{equation}
However (\ref{f(x)nu}) and (\ref{numbers}) yield that
if $z<x\leq z+\frac{3\Delta}{10}$, then
\begin{equation}
\label{intfsmall}
\int_{z}^{x}f<\frac{\nu}4.
\end{equation}
Since $T(z)\leq z$, (\ref{intgbig}) and (\ref{intfsmall})
yield that
if $z+\frac{2\Delta}{10}\leq x\leq z+\frac{3\Delta}{10}$, then
 $T(x)\leq z+\frac{\Delta}{10}$.
In particular
$$
\int_{z+\frac{2\Delta}{10}}^{z+\frac{3\Delta}{10}}(T(x)-x)^2f(x)\,dx
\geq \int_{z+\frac{2\Delta}{10}}^{z+\frac{3\Delta}{10}}
\left(\frac{\Delta}{10}\right)^2\frac{f(z)}2\,dx\gg \Delta^3f(z).
$$

\noindent{\bf Case 2 } $f(z)<g(z)\leq e^3 f(z)$.\\
Let $\lambda=(\frac{f(z)}{g(z)})^{\frac13}\geq 1/e$.
Since $2g(z)\leq 2e^3f(z)<50f(z)$ and
$\Delta=\frac{\nu\ln 2}{3f(z)}\,\ln\frac{g(z)}{f(z)}$,
if $z\leq x\leq z+\frac1{50}\Delta$, then
(\ref{f(x)nu}) and (\ref{g(x)nu}) yield
$$
\lambda\cdot f(z)\leq f(x)\leq \lambda^{-1}\cdot f(z)
\mbox{ \ and \ }
\lambda\cdot g(z)\leq g(x)\leq \lambda^{-1}\cdot g(z).
$$
 In particular
 if $z\leq s,t\leq z+\frac1{50}\Delta$, then
$\frac{f(s)}{g(t)}\leq\lambda$.
We deduce that if $z< x\leq z+\frac1{150}\Delta$ then
$$
\int_z^x f\leq  \int_z^{z+\lambda (x-z)} g.
$$
Thus $T(x)\leq z+\lambda (x-z)$ by $T(z)\leq z$, and hence
$$
x-T(x)\geq (1-\lambda) (x-z)= \lambda\left(\frac1{\lambda}-1\right)
(x-z)\geq
\frac{x-z}{3e}\cdot\ln\frac{g(z)}{f(z)}.
$$
It follows that
$$
\int_{z+\frac{\Delta}{300}}^{z+\frac{\Delta}{150}}(T(x)-x)^2f(x)\,dx
\gg \Delta^3f(z)\ln\frac{g(z)}{f(z)}.
$$

\noindent{\bf Case 3 } $e^{-3}f(z)\leq g(z)< f(z)$.\\
 Let
$\lambda=(\frac{f(z)}{g(z)})^{\frac13}\leq e$.
Since $\Delta=\frac{\nu\ln 2}{3f(z)}\,\ln\frac{f(z)}{g(z)}$,
if $z-\frac1{2}\,\Delta\leq x\leq z$, then
(\ref{f(x)nu}) and (\ref{g(x)nu}) yield
$$
\lambda^{-1}\cdot f(z)\leq f(x)\leq \lambda\cdot f(z)
\mbox{ \ and \ }
\lambda^{-1}\cdot g(z)\leq g(x)\leq \lambda\cdot g(z).
$$
 In particular
 if $z-\frac1{2}\Delta\leq s,t\leq z$, then
$\frac{f(s)}{g(t)}\geq \lambda$.
We deduce that if $z-\frac1{2e}\,\Delta< x\leq z$ then
$$
\int_x^z f\geq  \int_{z-\lambda (z-x)}^z g.
$$
Thus $T(x)\leq z-\lambda (z-x)$ by $T(z)\leq z$, and hence
$$
x-T(x)\geq (\lambda-1) (z-x)\geq
\frac{z-x}{3}\cdot\ln\frac{f(z)}{g(z)}.
$$
It follows that
$$
\int_{z-\frac{\Delta}{150}}^{z-\frac{\Delta}{300}}(T(x)-x)^2f(x)\,dx
\gg \Delta^3f(z)\ln\frac{f(z)}{g(z)}.
$$

\noindent{\bf Case 4 } $g(z)\leq e^{-3}f(z)$.\\
Since $\Delta=\frac{\nu\ln 2}{f(z)}$,
if $z-\Delta\leq x\leq z$, then
(\ref{f(x)nu}) and (\ref{g(x)nu}) yield that
$f(x)\geq  f(z)/2$ and $g(x)\leq 2g(z)$,
respectively.
 In particular
 if $z-\Delta\leq s,t\leq z$, then
$f(s)\geq 2g(t)$.
We deduce that if $z-\frac1{2}\,\Delta< x\leq z$ then
$$
\int_x^z f\geq  \int_{z-2(z-x)}^z g.
$$
Thus $T(x)\leq z-2(z-x)$ by $T(z)\leq z$, and hence
$x-T(x)\geq z-x$.
It follows that
$$
\int_{z-\frac{\Delta}{2}}^{z-\frac{\Delta}{4}}(T(x)-x)^2f(x)\,dx
\gg \Delta^3f(z).
\mbox{ \ }\proofbox
$$

\begin{theo}
\label{transL1}
If $f$ and $g$ are log-concave probability distributions on $\R$, and
$\int_{I_f}f(x)(T(x)-x)^2dx=\varepsilon\cdot f(w_f)^{-2}$
for $\varepsilon\in(0,1)$, then
$$
\int_\R|f(x)-g(x)|\,dx \ll
\sqrt[3]{\varepsilon}|\ln \varepsilon|^{\frac23}.
$$
\end{theo}
{\bf Remark } According to Example~\ref{exa3},
the exponent $\frac13$ of $\varepsilon$ is optimal
 in Lemma~\ref{transL1}.\\
\proof It is enough to
prove the statement if $\varepsilon<\varepsilon_0$,
where $\varepsilon_0\in(0,\frac12)$ is an absolute
constant specified later. We may assume that
$f(w_f)=1$, and hence
$f(x)\leq 2$ for any $x\in\R$ by Proposition~\ref{hwestimate} (iii),
and for the inverse $S$ of $T$,
$$
\int_{I_f}f(x)(T(x)-x)^2dx=
\int_{I_g}g(y)(S(y)-y)^2dy\leq \varepsilon.
$$
For $x\in\R$, we define
\begin{eqnarray*}
\nu(x)&=&\min\left\{\int_{-\infty}^{x}f,\int_x^{\infty}f\right\},\\
\tilde{\nu}(x)&=&\min\left\{\int_{-\infty}^{x}g,\int_x^{\infty}g\right\}.
\end{eqnarray*}

First we estimate $g$.
Since $\nu(w_f)=\frac12$, if $\varepsilon_0$ is small enough then
Propositions~\ref{transL11} and
\ref{transL10} yield that $\tilde{\nu}(w_f)> \frac14$ and
$g(w_f)\leq 2$, respectively. We conclude by
 Proposition~\ref{hwestimate} (ii) that $g(w_g)\leq 4$,
and hence $g(x)\leq 8$ for any $x\in\R$ by Proposition~\ref{hwestimate} (iii).

It follows by $f(x)\leq 2$ and Proposition~\ref{transL11}
that there exists a positive constant $c_1$
such that if $\nu(x)\geq c_1\sqrt[3]{\varepsilon}$
then $\tilde{\nu}(x)\geq \nu(x)/2$. Now applying Proposition~\ref{transL11}
to $g$,
and possibly increasing $c_1$, we have the following:
If $\nu(x)\geq c_1\sqrt[3]{\varepsilon}$
then $\tilde{\nu}(x)\leq 2\nu(x)$.
Finally, possibly increasing $c_1$ further,
if $\nu(x)\geq c_1\sqrt[3]{\varepsilon}$, then
$|\ln\frac{g(x)}{f(x)}|\leq \ln 2$ by Proposition~\ref{transL10}.
We choose $\varepsilon_0$ small enough
to satisfy $2c_1\sqrt[3]{\varepsilon_0}<\frac12$.

For $z\in\R$, we define
$\Delta(z)=\frac{\nu\ln 2}{450f(z)}\cdot |\ln\frac{g(z)}{f(z)}|$.
We assume $\nu(z)\geq c_1\sqrt[3]{\varepsilon}$,
and hence $\frac12\leq\frac{g(z)}{f(z)}\leq 2$.
It follows by Corollary~\ref{hwest} (i) that and $f(x)\geq f(z)/2$
$\nu(x)\leq 2\nu(z)$ if $x\in[z-\Delta(z),z+\Delta(z)]$.
We deduce using Proposition~\ref{transL10} and its
remark that there exists an absolute constant $c_2$
such that assuming $g(z)\neq f(z)$, we have
\begin{equation}
\label{shortinterval}
\int_{z-\Delta(z)}^{z+\Delta(z)}\frac{\nu(x)^2}{f(x)}\cdot
\left|\ln\frac{g(x)}{f(x)}\right|^3\,dx\leq c_2
\int_{z-\Delta(z)}^{z+\Delta(z)}(T(x)-x)^2f(x)\,dx.
\end{equation}

We define $z_1<z_2$ by the properties
$\nu(z_1)=\nu(z_2)=2c_1\sqrt[3]{\varepsilon}$.
We observe that if $g(z)\neq f(z)$
and some $x\in[z-\Delta(z),z+\Delta(z)]$
satisfies $\nu(x)\geq 2c_1\sqrt[3]{\varepsilon}$
then $\nu(z)\geq c_1\sqrt[3]{\varepsilon}$.
It is not hard to show based on (\ref{shortinterval})
that
$$
\int_{z_1}^{z_2}\frac{\nu(x)^2}{f(x)}\cdot
\left|\ln\frac{g(x)}{f(x)}\right|^3\,dx\leq c_2
\int_{\R}(T(x)-x)^2f(x)\,dx.
$$
Since $f(x)\leq 2$ and $\frac{|f(x)-g(x)|}{f(x)}\leq 4|\ln\frac{g(x)}{f(x)}|$
for $x\in [z_1,z_2]$, we deduce
\begin{eqnarray*}
\int_{z_1}^{z_2}\frac{\nu(x)^2|f(x)-g(x)|^3}{f(x)^2}\,dx&
= & 4\int_{z_1}^{z_2}\frac{\nu(x)^2}{f(x)}
\left(\frac{|f(x)-g(x)|}{f(x)}\right)^3\,dx\\
 & \leq &  4^4\int_{z_1}^{z_2}\frac{\nu(x)^2}{f(x)}
\left|\ln\frac{g(x)}{f(x)}\right|^3\,dx\leq 4^4c_2\varepsilon.
\end{eqnarray*}
It follows by the H\"older inequality  that
\begin{eqnarray*}
\int_{z_1}^{z_2}|f(x)-g(x)|\,dx&=&
\int_{z_1}^{z_2}\frac{\nu(x)^{\frac23}|f(x)-g(x)|}{f(x)^{\frac23}}
\cdot \frac{f(x)^{\frac23}}{\nu(x)^{\frac23}}\,dx\\
& \leq &
\left[\int_{z_1}^{z_2}\frac{\nu(x)^2|f(x)-g(x)|^3}{f(x)^2}\,dx\right]^{\frac13}
\times\\
&& \times\left[\int_{z_1}^{z_2}\frac{\nu(x)}{f(x)}\,dx\right]^{\frac23}.
\end{eqnarray*}
Here $f(x)=|\nu'(x)|$, therefore
\begin{eqnarray*}
\int_{z_1}^{z_2}|f(x)-g(x)|\,dx&\leq& (4^4c_2\varepsilon)^{\frac13}
\left[\int_{z_1}^{w_f}\frac{\nu'(x)}{\nu(x)}\,dx
+\int_{w_f}^{z_2}\frac{-\nu'(x)}{\nu(x)}\,dx\right]^{\frac23}\\
&=&(4^4c_2\varepsilon)^{\frac13}
\left[2\cdot \ln\mbox{$\frac12$}-2\cdot\ln(2c_1\sqrt[3]{\varepsilon})\right]^{\frac23}
\ll \sqrt[3]{\varepsilon}|\ln\varepsilon|^{\frac23}.
\end{eqnarray*}
On the other hand,
 $\tilde{\nu}(x_i)\leq 2\nu(x_i)=4c_1\sqrt[3]{\varepsilon}$, $i=1,2$,
yields that
$$
\int_{-\infty}^{z_1}|f(x)-g(x)|\,dx\leq 6c_1 \sqrt[3]{\varepsilon}
\mbox{ \ and \ }
\int_{z_2}^{\infty}|f(x)-g(x)|\,dx\leq 6c_1 \sqrt[3]{\varepsilon},
$$
and in turn we conclude Theorem~\ref{transL1}.
\proofbox

\section{The proof of Theorem~\ref{PLstab}}
\label{secmainproof}

For a non-negative, bounded, and not identically zero function $h$ on $\R$,
its log-concave hull is
$$
\tilde{h}(x)=\inf\{p(x):\mbox{ $p$ is a log-concave
function s.t. $h(t)\leq p(t)$ for $t\in\R$}\}.
$$
This $\tilde{h}$ is log-concave and $h(t)\leq \tilde{h}(t)$ for all $t\in\R$,
therefore we may take minimum in the definition.
Next we present a definition of $\tilde{h}$
in terms of $\ln h$. Let $J_h$ be the
set of all $x\in \R$ with $h(x)>0$, and let
$$
C_h=\{(x,y)\in\R^2:\,x\in J_h\mbox{ \ and \ }
y\leq \ln h(x)\}.
$$
This $C_h$ is convex if and only if $h$ is
log-concave. In addition
$J_{\tilde{h}}$ is the convex hull of $J_h$,
and  the interior of $C_{\tilde{h}}$
is the interior of the convex hull of $C_h$.
We also observe that for any unit vector $u\in\R^2$,
we have
\begin{equation}
\label{support}
\sup\{\langle u,v\rangle:\,v\in C_h\}=
\sup\{\langle u,v\rangle:\,v\in C_{\tilde{h}}\}.
\end{equation}

Let $f$, $g$ and $m$ be the functions in Theorem~\ref{PLstab}.
The condition of the Pr\'ekopa-Leindler inequality is equivalent with
\begin{equation}
\label{Ch}
\mbox{$\frac12$}(C_f+C_g)\subset C_m,
\end{equation}
where $C_f+C_g$ is
the Minkowski sum of the two sets.
Choose $x_0,y_0\in\R$ such that $f(x_0)>0$ and $g(y_0)>0$.
For any $x\in\R$, $m(\frac{x+x_0}2)\geq \sqrt{f(x_0)g(x)}$
and $m(\frac{x+y_0}2)\geq \sqrt{f(x)g(y_0)}$, and hence
$$
f(x)\leq \frac{m(\frac{x+y_0}2)^2}{g(y_0)}
\mbox{ \ and \ }
g(x)\leq \frac{m(\frac{x_0+x}2)^2}{f(x_0)}.
$$
Since $m$ is log-concave function with finite integral, it is bounded,
thus $f$ and $g$ are bounded, as well.
Therefore we may define the log-concave
hull of $f$ and $g$ of
$\tilde{f}$ and $\tilde{g}$, respectively.
It follows that $\tilde{f}(x)\geq f(x)$ and
$\tilde{g}(y)\geq g(y)$.
Since $m$ is log-concave, (\ref{support}) and (\ref{Ch})
yield that $m(\frac{x+y}2)\geq \sqrt{\tilde{f}(x)\tilde{g}(y)}$
for $x,y\in\R$.
We may assume that $\tilde{f}$ and $\tilde{g}$
are probability distributions
with zero mean, and $\tilde{f}(w_{\tilde{f}})=1$.
It follows that
\begin{equation}
\label{fgmest}
\int_\R f\geq 1-\varepsilon,
\mbox{ }  \int_\R g\geq 1-\varepsilon,
\mbox{ }\int_\R m\leq 1+\varepsilon.
\end{equation}

Next applying (\ref{stabcond}), Lemma~\ref{transdist}
and Theorem~\ref{transL1} to $\tilde{f}$ and $\tilde{g}$,
we conclude
\begin{equation}
\label{fgtilde}
\int_{\R}|\tilde{f}(t)-\tilde{g}(t)|\,dt\ll
\sqrt[3]{\varepsilon}|\ln \varepsilon|^{\frac43}.
\end{equation}
In addition (\ref{fgmest}) yields
\begin{equation}
\label{ffggtilde}
\int_{\R}|\tilde{f}(t)-f(t)|\,dt\leq
\varepsilon
\mbox{ \ and }
\int_{\R}|\tilde{g}(t)-g(t)|\,dt\leq
\varepsilon.
\end{equation}

Therefore to complete the proof of Theorem~\ref{PLstab},
all we have to do is to estimate $\int_{\R}|m(t)-\tilde{g}(t)|\,dt$.
For this, let $T:\,I_{\tilde{f}}\to I_{\tilde{g}}$ be
the transportation map satisfying
$$
\int_{-\infty}^x \tilde{f}(t)\,dt= \int_{-\infty}^{T(x)}\tilde{g}(t)\,dt.
$$
We note that $R(x)=\frac{x+T(x)}2$ is an increasing
and bijective map from $I_{\tilde{f}}$
into $\frac12(I_{\tilde{f}}+I_{\tilde{g}})$.
We define the function $h:\,\R\to\R$ as follows.
If $x\not\in\frac12(I_{\tilde{f}}+I_{\tilde{g}})$,
then $h(x)=0$,
and if $x\in I_{\tilde{f}}$, then
$$
h\left(\frac{x+T(x)}2\right)=
\sqrt{\tilde{f}(x)\tilde{g}(T(x))}.
$$
We have $h(x)\leq m(x)$,
and the proof of the Pr\'ekopa-Leindler inequality
using the transportation map in Section~\ref{sectrans}
shows that $\int_{\R}h\geq 1$.
We deduce by (\ref{fgmest}) that
\begin{equation}
\label{mh}
\int_{\R}|m(t)-h(t)|\,dt\leq
\varepsilon.
\end{equation}
To compare $h$ to $\tilde{g}$, we note that
$\int_{\R}h\leq 1+\varepsilon$ and (\ref{fgmest})
imply
\begin{equation}
\label{gh}
\int_{\R}h(t)-\tilde{g}(t)\,dt\leq 2\varepsilon.
\end{equation}
Let $B\subset \R$ be the set of all $t\in\R$
where $\tilde{g}(t)<h(t)$, and hence
$B\subset \frac12(I_{\tilde{f}}+I_{\tilde{g}})$.
In addition let $A=R^{-1}B\subset I_{\tilde{f}}$.
If $t=\frac{x+T(x)}2\in B$ for $x\in A$ then
as $\tilde{g}$ is concave and
$\tilde{f}(x)=\tilde{g}(T(x))T'(x)$, we have
\begin{eqnarray*}
[h(R(x)-\tilde{g}(R(x))]\cdot R'(x)&\leq &
\left[\sqrt{\tilde{f}(x)\tilde{g}(T(x))}-
\sqrt{\tilde{g}(x)\tilde{g}(T(x))}\right]\cdot\frac{1+T'(x)}2\\
&\leq& (\tilde{f}(x)-\tilde{g}(x))\cdot \frac{\sqrt{\tilde{g}(T(x))}}{\sqrt{\tilde{f}(x)}}
\cdot\frac{1+T'(x)}2\\
& =& (\tilde{f}(x)-\tilde{g}(x))\cdot
\left(1+\frac{(1-\sqrt{T'(x)})^2}{2\sqrt{T'(x)}}\right).
\end{eqnarray*}
In particular $\tilde{g}(x)<\tilde{f}(x)$ for $x\in A$.
It follows from (\ref{stabcond}) and (\ref{fgtilde}) that
\begin{eqnarray*}
\int_Bh(t)-\tilde{g}(t)\,dt&=&\int_A[h(R(x)-\tilde{g}(R(x))]\cdot R'(x)\,dx\\
&\leq&
\int_{I_{\tilde{f}(x)}}|\tilde{f}(x)-\tilde{g}(x)|+\tilde{f}(x)\cdot
\frac{(1-\sqrt{T'(x)})^2}{2\sqrt{T'(x)}}\,dx\\
&\ll& \sqrt[3]{\varepsilon}|\ln \varepsilon|^{\frac43}.
\end{eqnarray*}
It follows from (\ref{gh}) that
$\int_{\R}|h(t)-\tilde{g}(t)|\,dt\ll
\sqrt[3]{\varepsilon}|\ln \varepsilon|^{\frac43}$.
Therefore combining this estimate with (\ref{ffggtilde}) and (\ref{mh})
leads to $\int_{\R}|m(t)-g(t)|\,dt\ll
\sqrt[3]{\varepsilon}|\ln \varepsilon|^{\frac43}$.
In turn we deduce
$\int_{\R}|m(t)-f(t)|\,dt\ll
\sqrt[3]{\varepsilon}|\ln \varepsilon|^{\frac43}$
by (\ref{fgtilde}) and (\ref{ffggtilde}). \proofbox

\begin{remark}
Careful check of the argument
shows that the estimate for $\int_{\R}|m(t)-f(t)|\,dt$
and $\int_{\R}|m(t)-g(t)|\,dt$ is of the same order
as the estimate for $\int_{\R}|\tilde{f}(t)-\tilde{g}(t)|\,dt$.
Therefore to improve on the estimate in
Theorem~\ref{PLstab}, all one needs to improve
is (\ref{fgtilde}).
\end{remark}

\section{Appendix - Examples }

\begin{example}
\label{exa2}
If $f$ is an even log-concave probability distribution,
$g(x)=(1+\varepsilon)\cdot f((1+\varepsilon)x)$, and
$m(x)=(1+\varepsilon)\cdot f(x)$, then we have (\ref{stabcond}),
and
$$
\int_{I_f}f(x)(T(x)-x)^2dx=
\frac{\varepsilon^2}{(1+\varepsilon)^2}\int_{\R}x^2f(x)dx.
$$
\end{example}

\begin{example}
\label{exa3}
Let $f$ be the constant one on $[-\frac12,\frac12]$,
and let $g$ a modification such that
if $|x|\geq \frac12-\varepsilon$ then
$$
g(x)=e^{-\frac{|x|-\frac12+\varepsilon}{\varepsilon}}.
$$
In addition
$$
m(x)=\left\{
\begin{array}{rll}
1&\mbox{ }&  \mbox{ if $x\in [-\frac12,\frac12]$}  \\
e^{-\frac{|x|-\frac12}{\varepsilon}} &&\mbox{ otherwise}.
\end{array} \right.
$$
In this case $\int_{\R}m=1+\varepsilon$,
$$
\int_{\R}  f(x) \cdot\frac{(1-\sqrt{T'(x)})^2}{2\sqrt{T'(x)}}\,dx
\approx\varepsilon\mbox{ \ and \ }
\int_{\R} |f(x)-g(x)|\,dx\approx\varepsilon.
$$
Moreover $\int_{\R}f(x)(T'(x)-1)^2dx=\infty$ and
$\int_{\R}f(x)(T(x)-x)^2dx\approx \varepsilon^3$.
\end{example}

\noindent{\bf Acknowledgement: } We are grateful for the help of
Marianna Cs\"ornyei, Katalin Marton, Franck Barthe, Francesco
Maggi, Cedric Villani, in the preparation of this manuscript.

\end{document}